\documentclass[letterpaper, 10 pt, conference]{ieeeconf}  

\IEEEoverridecommandlockouts                              
\usepackage{graphicx,amssymb,amstext}
\usepackage{amsmath}
\usepackage{tikz,lscape}
\usepackage{hyperref}
\usepackage{epsfig}
\usetikzlibrary{shapes,arrows}
\usepackage{verbatim}
\usepackage{cite}

\newtheorem{theorem}{Theorem}

\title{\LARGE \bf
Robust Smoothing for Estimating Optical Phase\\
Varying as a Continuous Resonant Process
}

\author{Shibdas Roy*, Obaid Ur Rehman, Ian R. Petersen and Elanor H. Huntington%
\thanks{This work was supported by the Australian Research Council (ARC).}%
\thanks{The authors are with the School of Engineering and Information Technology, University of New South Wales at the Australian Defence Force Academy, Canberra, ACT 2600, Australia.}%
\thanks{*\tt\small shibdas.roy at student.adfa.edu.au}%
}

\begin{document}


\maketitle
\thispagestyle{empty}
\pagestyle{empty}

\begin{abstract}

Continuous phase estimation is known to be superior in accuracy as compared to static estimation. The estimation process is, however, desired to be made robust to uncertainties in the underlying parameters. Here, homodyne phase estimation of coherent and squeezed states of light, evolving continuously under the influence of a second-order resonant noise process, are made robust to parameter uncertainties using a robust fixed-interval smoother, designed for uncertain systems satisfying a certain integral quadratic constraint. We observe that such a robust smoother provides improved worst-case performance over the optimal smoother and also performs better than a robust filter for the uncertain system.

\end{abstract}

\section{INTRODUCTION}

\bstctlcite{BSTcontrol}

Quantum phase estimation \cite{WM,GLM2} is central to various fields such as metrology \cite{GLM1}, quantum computation \cite{HWA}, communication \cite{SPK,CHD} and quantum cryptography \cite{IWY}. Real-time feedback in homodyne estimation of a \emph{static} unknown phase can yield mean-square errors much lower than without using feedback \cite{HMW,WK1,WK2,MA,BW1}. However, it is experimentally more relevant to precisely estimate a phase, that is \emph{continuously varying} under the influence of an unmeasured classical noise process \cite{BW2,TSL,TW,YNW}. A classical process coupled dynamically to a quantum system under continuous measurement may be estimated in various ways: prediction, filtering or smoothing \cite{MT}. In particular, smoothing uses both past and future measurements to yield a more accurate estimate than filtering alone, that uses only past measurements.

The \emph{fixed-interval} smoothing problem \cite{LK,JSM,WWS} considers measurements over a fixed time-interval $\tau$, where the estimation time is $t=\tau-q$ for some $q:0<q<\tau$. One solution, the Mayne-Fraser two-filter smoother \cite{DQM,DCF,RKM}, consists of a forward-time Kalman filter and also a backward-time Kalman filter called an ``information filter'' \cite{FP} and it combines the two estimates to yield the optimal smoothed estimate. The information filter and the smoother were combined into a single backward smoother by Rauch, Tung and Striebel (RTS) \cite{RTS}.

In \cite{TW,YNW}, the signal phase to be estimated is allowed to evolve under the influence of an unmeasured continuous-time Ornstein-Uhlenbeck (OU) noise process. However, the estimation process heavily relies on the underlying parameters being precisely known, which is practically not feasible due to unavoidable external noises and/or apparatus imperfections. It is, therefore, desired to make the estimation process robust to uncertainties in these parameters to improve the quality of the phase measurement. In \cite{RPH1} the authors have shown that a robust guaranteed cost filter \cite{PM} yields better worst-case accuracy in the phase estimate than an optimal filter with parameter uncertainty in the phase being measured. The authors have also demonstrated improvement in continuous phase estimation with robust fixed-interval smoothing \cite{MSP} for coherent \cite{RPH2} and squeezed \cite{RPH3} states.

These works, however, considered a simplistic OU noise process modulating the signal phase to be estimated, and the kind of noises that in practice corrupt the signal are more complicated than an OU process. The authors have illustrated in \cite{RPH4} a guaranteed cost robust filter designed for a more complicated and practically relevant second-order resonant noise process. The resonant noise process is generated by a piezo-electric transducer with a linear transfer function, driven by a white noise, e.g. see \cite{IMY}. Here, we aim to extend robustness for such a resonant noise process to consider smoothing that is known to provide with more precise estimates than filtering alone. Such a robust smoother is essentially acausal and used for offline estimation, as mentioned before, and is useful for highly sensitive measurements requiring better precision than real-time estimation, as in the case of gravitational wave detection. This is in contrast to \cite{RPH4}, where only the feedback filter was made robust and, therefore, yielded only partially better estimates as a robust smoother with the same uncertainty in the system, although the robust filter can be used in real-time.

We design a robust fixed-interval smoother \cite{MSP}, meant for uncertain systems admitting a certain integral quadratic constraint, for estimation of the phase of a coherent light beam phase-modulated by such a resonant noise process. We compare the performance of the robust and RTS smoothers with the robust and Kalman filters for the uncertain system. We also extend this by designing such a robust fixed-interval smoother for estimating the phase of a phase-squeezed light beam and compare its behaviour with a corresponding RTS smoother as well as the robust and Kalman filters for the uncertain system. Note that this work is significantly different from \cite{RPH2} and \cite{RPH3} since here we consider a more realistic noise process, modulating the phase to be estimated, and in doing so, we get better improvement in the performance and sensitivity of the robust smoother with system uncertainty.

\section{ROBUST FIXED-INTERVAL SMOOTHING}\label{sec:unc}

In this section, we outline the robust fixed-interval smoothing theory from \cite{MSP}, that we use later in this paper. Consider an uncertain system described by the state equations

\small
\begin{equation}\label{eq:unc_0}
\begin{split}
\dot{x}(t) &= [A(t) + B_1(t)\Delta_1(t)K(t)]x(t) + B_1w(t) + B_2(t)u(t),\\
y(t) &= [C(t) + \Delta_2(t)K(t)]x(t) + v(t),\\
z(t) &= K(t)x(t) + G(t)u(t),
\end{split}
\end{equation}\normalsize
where $x(t)$ is the state, $y(t)$ is the measured output, $u(t)$ is a known input, $z(t)$ is the uncertainty output, $w(t)$ and $v(t)$ are white noises. $A(\cdot), B(\cdot), K(\cdot)$ and $C(\cdot)$ are bounded piecewise continuous matrix functions. Furthermore, $\Delta_1(t)$ and $\Delta_2(t)$ are uncertainty matrices satisfying
\begin{equation}\label{eq:unc_1}
||\left[\begin{array}{cc} \Delta_1(t)^TQ(t)^{\frac{1}{2}} & \Delta_2(t)^TR(t)^{\frac{1}{2}} \end{array}\right]|| \leq 1
\end{equation}
for all $t$, where $Q(\cdot) = Q(\cdot)^T$ and $R(\cdot) = R(\cdot)^T$ are bounded piecewise continuous matrix functions, satisfying $Q(t) \geq \delta I$, $R(t) \geq \delta I$ for all $t$, for some constant $\delta > 0$.

Let $X_0 = X_0^T > 0$ be a given matrix, $x_0$ be a given real vector, $d > 0$ be a given constant. Then, we require the initial conditions $x(0)$ to satisfy the inequality
\begin{equation}\label{eq:unc_2}
(x(0)-x_0)^TX_0(x(0)-x_0) \leq d.
\end{equation}

Moreover, $\tilde{w}(t)$ and $\tilde{v}(t)$ are uncertainty inputs, given by
\begin{equation}\label{eq:unc_3}
\begin{split}
\tilde{w}(t) &= \Delta_1(t)[K(t)x(t)+G(t)u(t)] + w(t),\\
\tilde{v}(t) &= \Delta_2(t)[K(t)x(t)+G(t)u(t)] + v(t)
\end{split}
\end{equation}

Then, for a given finite time interval $[0, \tau]$, (\ref{eq:unc_1}), (\ref{eq:unc_2}), (\ref{eq:unc_3}) constitute the following integral quadratic constraint (IQC) as the description of uncertainty for the system (\ref{eq:unc_0}):
\vspace*{-2mm}
\begin{equation}\label{eq:unc_iqc}
\begin{split}
(x(0)&-x_0)^TX_0(x(0)-x_0) + \int_0^\tau (\tilde{w}(t)^TQ(t)\tilde{w}(t)\\
&+\tilde{v}(t)^TR(t)\tilde{v}(t))dt \leq d + \int_0^\tau ||z(t)||^2 dt.
\end{split}
\end{equation}

A solution to the robust fixed-interval smoothing problem for this uncertain system involves the Riccati equations:

\vspace*{-2mm}
\small
\begin{equation}\label{eq:unc_ric1}
\begin{split}
-\dot{X}(t) = &X(t)A(t) + A(t)^TX(t) + X(t)B_1(t)Q(t)^{-1}B_1(t)^TX(t)\\
            &+K(t)^TK(t) - C(t)^TR(t)C(t); \quad X(0) = X_0,
\end{split}
\end{equation}
\begin{equation}\label{eq:unc_ric2}
\begin{split}
-\dot{Y}(t) = &Y(t)A(t) + A(t)^TY(t) + Y(t)B_1(t)Q(t)^{-1}B_1(t)^TY(t)\\
            &-K(t)^TK(t) + C(t)^TR(t)C(t); \quad Y(\tau) = 0.
\end{split}
\end{equation}\normalsize

It will also include a solution to the differential equations:

\vspace*{-2mm}
\small
\begin{equation}\label{eq:unc_diff1}
\begin{split}
\dot{\eta}(t) &= -[A(t)+B_1(t)Q(t)^{-1}B_1(t)^TX(t)]^T\eta(t)+C(t)^TR(t)\\
			  &\times y_0(t)+[K(t)^TG(t)+X(t)B_2(t)]u_0(t); \quad \eta(0) = X_0x_0
\end{split}
\end{equation}\normalsize
for $t \in [0,\tau -q]$ and

\vspace*{-2mm}
\small
\begin{equation}\label{eq:unc_diff2}
\begin{split}
-\dot{\xi}(t) &= [A(t)-B_1(t)Q(t)^{-1}B_1(t)^TY(t)]^T\xi(t)+C(t)^TR(t)\\
			  &\times y_0(t)-[Y(t)B_2(t)-K(t)^TG(t)]u_0(t); \quad \xi(\tau) = 0
\end{split}
\end{equation}\normalsize
for $t \in [\tau -q,\tau]$.

\begin{theorem}
Assume that (\ref{eq:unc_ric1}) has a solution over time interval $t \in [0,\tau -q]$ such that $X(\tau -q) > 0$ and (\ref{eq:unc_ric2}) has a solution over time interval $t \in [\tau -q,\tau]$ such that $Y(\tau -q) > 0$. Then, the set $X_{\tau -q}[x_0,u_0(\cdot)|_0^\tau ,y_0(\cdot)|_0^\tau ,d]$ of all possible states $x(\tau -q)$ at time $\tau -q$ for the uncertain system (\ref{eq:unc_0}) with uncertainty inputs and initial conditions satisfying (\ref{eq:unc_iqc}) is bounded and is given by:

\vspace*{-2mm}
\small
\begin{equation}\label{eq:unc_ellipse}
\begin{split}
X_{\tau -q}&[x_0,u_0(\cdot)|_0^\tau ,d] = \left\lbrace x_{\tau -q}: x_{\tau -q}^TX(\tau -q)x_{\tau -q}\right. \\
&-2x_{\tau -q}^T\eta(\tau -q) + h_{\tau -q} + x_{\tau -q}^TY(\tau -q)x_{\tau -q}\\
&\left.-2x_{\tau -q}^T\xi(\tau -q) + s_{\tau -q} \leq d \right\rbrace
\end{split}
\end{equation}\normalsize
where $\eta(t)$ and $\xi(t)$ are solutions to (\ref{eq:unc_diff1}) and (\ref{eq:unc_diff2}) and
\vspace*{-2mm}
\small
\begin{equation}
\begin{split}
h_{\tau -q} &= x_0^TX_0x_0 + \int_0^{\tau -q}\left\{y_0(t)^TR(t)y_0(t) - u_0(t)^TG(t)^T\right.\\
&\times G(t)u_0(t) - \eta(t)^TB_1(t)Q(t)^{-1}B_1(t)^T\eta(t)\\
&\left.+2u_0(t)^TB_2(t)\eta(t)\right\}dt,
\end{split}
\end{equation}
\vspace*{-2mm}
\begin{equation}
\begin{split}
s_{\tau -q} &= \int_0^{\tau -q}\left\{y_0(t)^TR(t)y_0(t) - u_0(t)^TG(t)^TG(t)u_0(t)\right.\\
&\left.-\xi(t)^TB_1(t)Q(t)^{-1}B_1(t)^T\xi(t)-2u_0(t)^TB_2(t)\xi(t)\right\}dt
\end{split}
\end{equation}\normalsize
\end{theorem}

\section{RESONANT NOISE PROCESS}\label{sec:resonant}

The resonant noise process under consideration in this paper is typically generated by a piezo-electric transducer (PZT) driven by an input white noise. The simplified transfer function of a typical PZT is given by:
\begin{equation}\label{eq:pzt_tf}
{G(s) := \frac{\phi}{v} = \frac{\kappa}{s^2+2\zeta\omega_r s+\omega_r^2},}
\end{equation}
where $\kappa$ is the gain, $\zeta$ is the damping factor, $\omega_r$ is the resonant frequency (rad/s), $v$ is a zero-mean white Gaussian noise with unity amplitude and $\phi$ is the PZT output that modulates the phase to be estimated.

\begin{figure}[!b]
\hspace*{-5mm}
\includegraphics[width=0.56\textwidth]{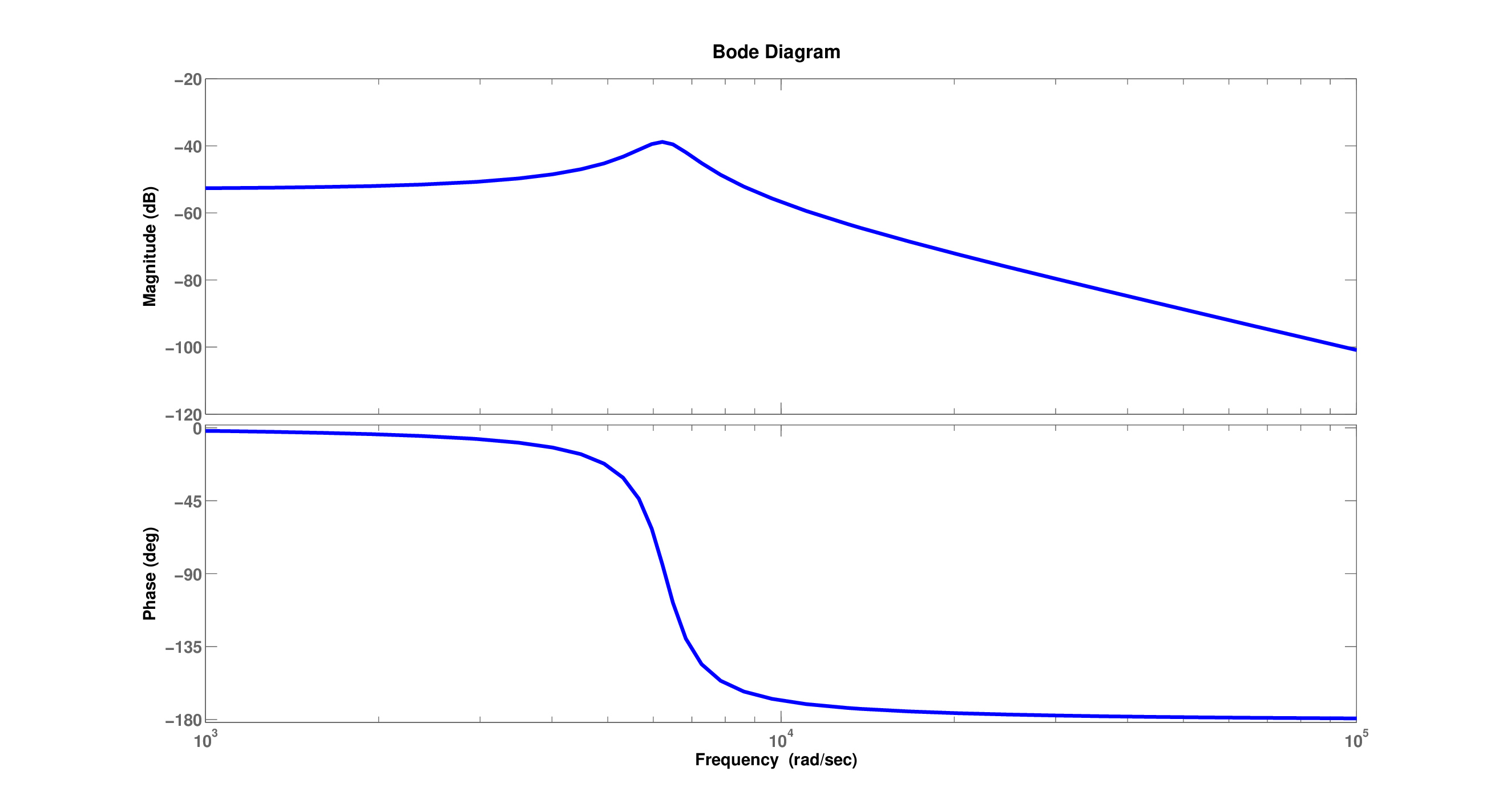}
\caption{Bode plot of the resonant noise process transfer function.}
\label{fig:bode_resonant}
\end{figure}

We use the following values for the parameters above: $\kappa = 9 \times 10^4$, $\zeta = 0.1$ and $\omega_r = 6.283 \times 10^3$ rad/s (1 kHz). Fig. \ref{fig:bode_resonant} shows the Bode plot of the transfer function (\ref{eq:pzt_tf}).

Let $x_1 := \phi$ and $x_2 := \dot{\phi}$. A state-space realization of the transfer function (\ref{eq:pzt_tf}) is:
\begin{equation}\label{eq:process_eqn}
{\mathbf{\dot{x}} = \mathbf{Ax} + \mathbf{G}v,}
\end{equation}
where
\[
\mathbf{x} := \left[\begin{array}{c}
x_1\\
x_2
\end{array}\right], \quad
\mathbf{A} := \left[\begin{array}{cc}
0 & 1\\
-\omega_r^2 & -2\zeta\omega_r
\end{array}\right], \quad
\mathbf{G} := \left[\begin{array}{c}
0\\
\kappa
\end{array}\right].
\]

\section{OPTIMAL SMOOTHER}\label{sec:optimal_smoother}
The homodyne ``smoothed'' phase estimation of a weak coherent state of light is optimal for given values of the parameters when using an offline RTS estimator alongwith a Kalman filter in the feedback loop to adjust the phase of the local oscillator. Under a linearization approximation, the output homodyne photocurrent $I(t)$ is \cite{TW}:
\begin{equation}
I(t)dt = 2|\alpha|[\phi(t)-\hat{\phi}(t)]dt + dW(t),
\end{equation}
where $|\alpha|$ is the amplitude of the coherent state, $\hat{\phi}$ is the \emph{intermediate phase estimate}, and $W(t)$ is a Wiener process arising from the quantum vacuum fluctuations.

The \emph{instantaneous estimate} $\theta(t)$ is defined as \cite{TW}:
\begin{equation}
\theta(t) := 2|\alpha |\hat{\phi}(t) + I(t) = 2|\alpha |\phi(t) + \frac{dW(t)}{dt}.
\end{equation}
\begin{equation}
\therefore \qquad \theta = \mathbf{Hx}+\mathbf{J}w,
\end{equation}
where $w := \frac{dW}{dt}$ is zero-mean Gaussian white noise with unity amplitude, $\mathbf{H} := \left[\begin{array}{cc} 2|\alpha| & 0 \end{array}\right]$ and $\mathbf{J} := 1$. Thus,
\begin{eqnarray}\label{eq:sys_model}
{\begin{split}
\textsf{\small Process model:} \ \ \mathbf{\dot{x}} &= \mathbf{Ax} + \mathbf{G}v, \\
\textsf{\small Measurement model:} \ \ \theta &= \mathbf{Hx} + \mathbf{J}w,
\end{split}}
\end{eqnarray}
where $E[v(t)v(t_1)] = \mathbf{N}\delta(t - t_1),
E[w(t)w(t_1)] = \mathbf{S}\delta(t - t_1),
E[v(t)w(t_1)] = 0$. Since $v$ and $w$ are unity amplitude white noise processes, both $\mathbf{N}$ and $\mathbf{S}$ are unity (scalars).

\subsection{Forward Kalman Filter}

The continuous-time algebraic Riccati equation to solve to construct the steady-state forward Kalman filter is \cite{RGB}:
\begin{equation}\label{eq:kalman_riccati}
\mathbf{AP_f}+\mathbf{P_fA}^T+\mathbf{GNG}^T-\mathbf{P_fH}^T(\mathbf{JSJ}^T)^{-1}\mathbf{HP_f} = \mathbf{0},
\end{equation}
where $\mathbf{P_f}$ is the error-covariance of the forward filter.

The forward Kalman filter equation is \cite{RGB}:
\begin{equation}\label{eq:fwd_kalman_filter}
{\mathbf{\dot{\hat{x}}} = (\mathbf{A}-\mathbf{K_fH})\mathbf{\hat{x}}+\mathbf{K_fHx}+\mathbf{K_fJ}w,}
\end{equation}
where $\mathbf{K_f} := \mathbf{P_fH}^T(\mathbf{JSJ}^T)^{-1}$ is the Kalman gain.

\subsection{Backward Kalman Filter}
The continuous-time algebraic Riccati equation to solve to construct the steady-state backward Kalman filter is \cite{LXP}:

\vspace*{-2mm}
\small
\begin{equation}
\mathbf{-AP_b}-\mathbf{P_bA}^T+\mathbf{GNG}^T-\mathbf{P_bH}^T(\mathbf{JSJ}^T)^{-1}\mathbf{HP_b} = \mathbf{0},
\end{equation}
\normalsize
where $\mathbf{P_b}$ is the error-covariance of the backward filter.

The backward Kalman filter equation is \cite{LXP}:
\begin{equation}\label{eq:bwd_kalman_filter}
{\mathbf{\dot{\hat{x}}} = (-\mathbf{A}-\mathbf{K_bH})\mathbf{\hat{x}}+\mathbf{K_bHx}+\mathbf{K_bJ}w,}
\end{equation}
where $\mathbf{K_b} := \mathbf{P_bH}^T(\mathbf{JSJ}^T)^{-1}$ is the Kalman gain.

\subsection{Smoother Error}
The smoother error covariance matrix $\mathbf{P_s}$ is obtained as:
\begin{equation}\label{eq:kalman_smoother_error}
\mathbf{P_s} = \left(\mathbf{P_f}^{-1}+\mathbf{P_b}^{-1}\right)^{-1},
\end{equation}
since the forward and backward estimates of the optimal RTS smoother are independent \cite{LXP}.

Using $\kappa = 9 \times 10^4$, $\zeta = 0.1$ and $\omega_r = 6.283 \times 10^3$ rad/s (1 kHz) as in Section \ref{sec:resonant} and $|\alpha| = 5 \times 10^2$, we get:

\vspace*{-2mm}
\small
\begin{equation}
\mathbf{P_s} = \left[\begin{array}{cc}
3.7748607 \times 10^{-3} & -7.2880146 \times 10^{-15}\\
-7.2880146 \times 10^{-15} & 3.7098537 \times 10^5
\end{array}\right].
\end{equation}\normalsize

\section{ROBUST SMOOTHER}

In this section, we make our smoother robust to uncertainty in the resonant frequency $\omega_r$ underlying the system matrix $\mathbf{A}$ using robust fixed-interval smoothing approach from \cite{MSP}, as outlined in section \ref{sec:unc}.

We introduce uncertainty in $\mathbf{A}$ as follows:
\[
\mathbf{A} \to \mathbf{A} + \left[\begin{array}{cc}
0 & 0\\
-\mu\delta\omega_r^2 & 0
\end{array}\right],
\]
where uncertainty is introduced in the resonant frequency $\omega_r$ through $\delta$. Furthermore, $\mathbf{\Delta} := \left[\begin{array}{cc}\delta & 0 \end{array}\right]$ is an uncertain parameter satisfying $||\mathbf{\Delta}|| \leq 1$ which implies $\delta^2 \leq 1$. Moreover, $\mu \in [0,1)$ determines the level of uncertainty. Uncertainty in $\zeta$ is deliberately not included here, since we do not get significant improvement in estimation error with robust smoother over that with RTS smoother in this case.

The process and measurement models of (\ref{eq:sys_model}) become:
\begin{eqnarray}\label{eq:uncertain_model}
{
\begin{split}
\textsf{\small Process model:} \ \ \mathbf{\dot{x}} &= (\mathbf{A}+\mathbf{G\Delta K})\mathbf{x} + \mathbf{G}v, \\
\textsf{\small Measurement model:} \ \ \theta &= \mathbf{H}\mathbf{x} + \mathbf{J}w,
\end{split}}
\end{eqnarray}
where 
$\mathbf{K}:=\left[\begin{array}{cc}
-\frac{\mu\omega_r^2}{\kappa} & 0\\
0 & 0
\end{array}\right]$, so $\mathbf{G\Delta K} = \left[\begin{array}{cc}
0 & 0\\
-\mu\delta\omega_r^2 & 0
\end{array}\right]$.

The IQC of (\ref{eq:unc_iqc}) for our case is:
\begin{equation}
\int_0^\tau (\tilde{w}^2+\tilde{v}^2)dt \leq 1 + \int_0^\tau ||\mathbf{z}||^2dt,
\end{equation}
where $\mathbf{z} = \mathbf{Kx}$ is the \emph{uncertainty output}, and $\tilde{w} = \mathbf{\Delta Kx}+v$ and $\tilde{v}=w$ are the \emph{uncertainty inputs}. Here, $X_0 = 0$, since no \emph{a-priori} information exists about the initial condition of the state in our case. Also, $d = 1$, since the amplitudes of the white noise processes $v$ and $w$ have been assumed to be unity. Thus, we would have $\mathbf{Q}=1$ and $\mathbf{R}=1$ in our case.

The steady-state forward Riccati equation from (\ref{eq:unc_ric1}):
\begin{eqnarray}
\mathbf{YA}+\mathbf{A}^T\mathbf{Y}+\mathbf{YG}\mathbf{Q}^{-1}\mathbf{G}^T\mathbf{Y}
+\mathbf{K}^T\mathbf{K}-\mathbf{H}^T\mathbf{RH}=\mathbf{0},
\end{eqnarray}
where we have used $\mathbf{Y}$ in place of $X(t)$.

The steady-state backward Riccati equation from (\ref{eq:unc_ric2}):
\begin{eqnarray}
\mathbf{ZA}+\mathbf{A}^T\mathbf{Z}-\mathbf{ZG}\mathbf{Q}^{-1}\mathbf{G}^T\mathbf{Z}
-\mathbf{K}^T\mathbf{K}+\mathbf{H}^T\mathbf{RH}=\mathbf{0},
\end{eqnarray}
where we have used $\mathbf{Z}$ in place of $Y(t)$.

Next, (\ref{eq:unc_diff1}) in this case yields:
\begin{equation}
\mathbf{\dot{\eta}}=-(\mathbf{A}+\mathbf{GQ}^{-1}\mathbf{G}^T\mathbf{Y})^T\mathbf{\eta} + \mathbf{H}^T\mathbf{R}\theta.
\end{equation}

Likewise, (\ref{eq:unc_diff2}) for reverse-time in this case yields:
\begin{equation}
\mathbf{\dot{\xi}}=(\mathbf{A}-\mathbf{GQ}^{-1}\mathbf{G}^T\mathbf{Z})^T\mathbf{\xi} + \mathbf{H}^T\mathbf{R}\theta.
\end{equation}

The forward filter is, then, simply: $\mathbf{\hat{x}_f}=\mathbf{Y}^{-1}\mathbf{\eta}$. Likewise, the backward filter is: $\mathbf{\hat{x}_b}=\mathbf{Z}^{-1}\mathbf{\xi}$.

The robust smoother for the uncertain system would, then, be the centre of the ellipse of (\ref{eq:unc_ellipse}):
\begin{equation}\label{eq:robust_smoother_eqn}
\mathbf{\hat{x}} = (\mathbf{Y}+\mathbf{Z})^{-1}(\mathbf{\eta} - \mathbf{\xi}).
\end{equation}

\section{COMPARISON OF THE SMOOTHERS}

\subsection{Error Analysis}\label{sec:lyap_method}

\subsubsection{Forward Filter}

We augment the system given by (\ref{eq:uncertain_model}) with the forward Kalman filter (\ref{eq:fwd_kalman_filter}) and represent the augmented system by the state-space model:
\begin{equation}\label{eq:ss_model}
\mathbf{\dot{\overline{x}}} = \mathbf{\overline{A}\, \overline{x}} + \mathbf{\overline{B}\, \overline{w}},
\end{equation}
where

\qquad \qquad \( \mathbf{\overline{x}} := 
\left[ \begin{array}{c}
\mathbf{x} \\
\mathbf{\hat{x}}
\end{array} \right] \)
\qquad and \qquad
\( \mathbf{\overline{w}} := 
\left[ \begin{array}{c}
v \\
w
\end{array} \right]. \)

\vspace*{2mm}
\[\therefore \mathbf{\overline{A}} = 
\left[ \begin{array}{cc}
\mathbf{A}+\mathbf{G\Delta K} & \mathbf{0} \\
\mathbf{K_fH} & \mathbf{A}-\mathbf{K_fH}
\end{array} \right],
\mathbf{\overline{B}} = 
\left[ \begin{array}{cc}
\mathbf{G} & \mathbf{0} \\
\mathbf{0} & \mathbf{K_fJ}
\end{array} \right]. \]

\vspace*{2mm}

For the continuous-time state-space model (\ref{eq:ss_model}), the steady-state state covariance matrix $\mathbf{P_{fs}}$ is obtained by solving the Lyapunov equation:
\begin{equation}\label{eq:lyapunov}
\mathbf{\overline{A}P_{fs}} + \mathbf{P_{fs}}\mathbf{\overline{A}}^T + \mathbf{\overline{B}\, \overline{B}}^T = \mathbf{0},
\end{equation}
where $\mathbf{P_{fs}}$ is the symmetric matrix
\[ 
\mathbf{P_{fs}} := E(\mathbf{\overline{x}\, \overline{x}}^T) :=
\left[ \begin{array}{cc}
\mathbf{\Sigma} & \mathbf{M_f}\\
\mathbf{M_f}^T & \mathbf{N_f}
\end{array} \right].
\]

The state estimation error can be written as:
\[ \mathbf{e_f} := \mathbf{x} - \mathbf{\hat{x}} = [\mathbf{1} \, -\mathbf{1}]\mathbf{\overline{x}}, \]
which is mean zero since all of the quantities determining $\mathbf{e_f}$ are mean zero.

The error covariance matrix is then given as:
\begin{equation}
\begin{split}
\mathbf{\Pi_f} :&= E(\mathbf{e_fe_f}^T) = [\mathbf{1} \, -\mathbf{1}]E(\mathbf{\overline{x}\, \overline{x}}^T)
\left[ \begin{array}{c}
\mathbf{1} \\
-\mathbf{1}
\end{array} \right] \\ 
&=
[\mathbf{1} \, -\mathbf{1}]
\left[ \begin{array}{cc}
\mathbf{\Sigma} & \mathbf{M_f} \\
\mathbf{M_f}^T & \mathbf{N_f}
\end{array} \right]
\left[ \begin{array}{c}
\mathbf{1} \\
-\mathbf{1}
\end{array} \right] \\ 
&= \mathbf{\Sigma} - \mathbf{M_f} - \mathbf{M_f}^T + \mathbf{N_f}.
\end{split}
\end{equation}

Since we are mainly interested in estimating $x_1 = \phi$, the estimation error covariance of interest is $\mathbf{\Pi_f}(1,1)$.

\subsubsection{Backward Filter}

The augmented system state-space model (\ref{eq:ss_model}) for the backward Kalman filter (\ref{eq:bwd_kalman_filter}) would have:

\[ \mathbf{\overline{A}} = 
\left[ \begin{array}{cc}
\mathbf{A}+\mathbf{G\Delta K} & \mathbf{0} \\
\mathbf{K_bH} & -\mathbf{A}-\mathbf{K_bH}
\end{array} \right],
\mathbf{\overline{B}} = 
\left[ \begin{array}{cc}
\mathbf{G} & \mathbf{0} \\
\mathbf{0} & \mathbf{K_bJ}
\end{array} \right]. \]

\vspace*{2mm}
As pointed out in \cite{RPH2}, in the steady-state case, the reverse time output process is also a stationary random process with the same auto-correlation function as the forward time output process. So, it can be regarded as being generated by the same process that generated the forward time output process. This is why we augment the forward-time process equation with the backward time Kalman filter above.

We then solve (\ref{eq:lyapunov}), with $\mathbf{P_{fs}}$ replaced by
\[ 
\mathbf{P_{bs}} := E(\mathbf{\overline{x}\, \overline{x}}^T) :=
\left[ \begin{array}{cc}
\mathbf{\Sigma} & \mathbf{M_b} \\
\mathbf{M_b}^T & \mathbf{N_b}
\end{array} \right],
\]
for the backward filter.

The error covariance matrix is, thus:
\begin{equation}
\mathbf{\Pi_b} := E(\mathbf{e_be_b}^T) = \mathbf{\Sigma} - \mathbf{M_b} - \mathbf{M_b}^T + \mathbf{N_b}.
\end{equation}

Here, the error of interest is $\mathbf{\Pi_b}(1,1)$.

\subsubsection{Cross-Correlation Term}

The forward and backward estimates are not independent and are correlated in this case, unlike in Section \ref{sec:optimal_smoother}. The cross-correlation term is \cite{RPH2}:
\begin{equation}
\mathbf{\Pi_{fb}} := E(\mathbf{e_fe_b}^T) = \mathbf{\Sigma} - \mathbf{M_f}^T - \mathbf{M_b} + \mathbf{\alpha\Sigma\beta},
\end{equation}
where $\mathbf{\alpha} := \mathbf{M_f}^T\mathbf{\Sigma}^{-1}$ and $\mathbf{\beta} := \mathbf{\Sigma}^{-1}\mathbf{M_b}$ \cite{WWS}. Here, the error of interest is $\mathbf{\Pi_{fb}}(1,1)$.

\subsubsection{Smoother Error}

The smoother error covariance of interest $\sigma^2$ is \cite{RPH2}:
\begin{equation}
\Pi := \frac{\mathbf{\Pi_f}(1,1)\mathbf{\Pi_b}(1,1)-\mathbf{\Pi_{fb}}(1,1)^2}{\mathbf{\Pi_f}(1,1) + \mathbf{\Pi_b}(1,1) - 2\mathbf{\Pi_{fb}}(1,1)}.
\end{equation}

\subsection{Comparison of Estimation Errors}

\begin{figure}[!b]
\hspace*{-5mm}
\includegraphics[width=0.56\textwidth]{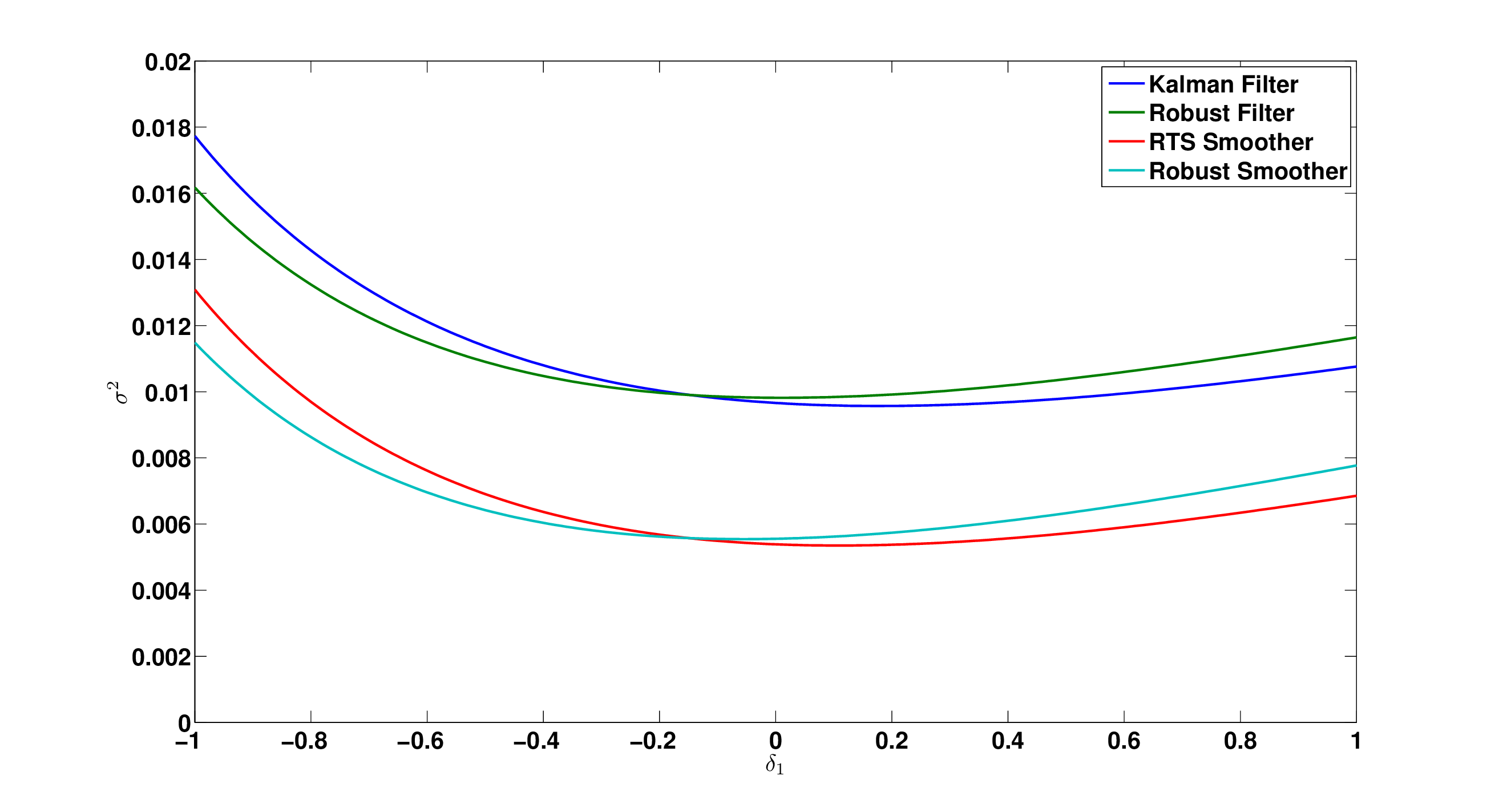}
\caption{Coherent State: Comparison of the smoothers for $\mu=0.5$.}
\label{fig:sql_mu_50}
\end{figure}

\begin{figure}[!htb]
\hspace*{-5mm}
\includegraphics[width=0.56\textwidth]{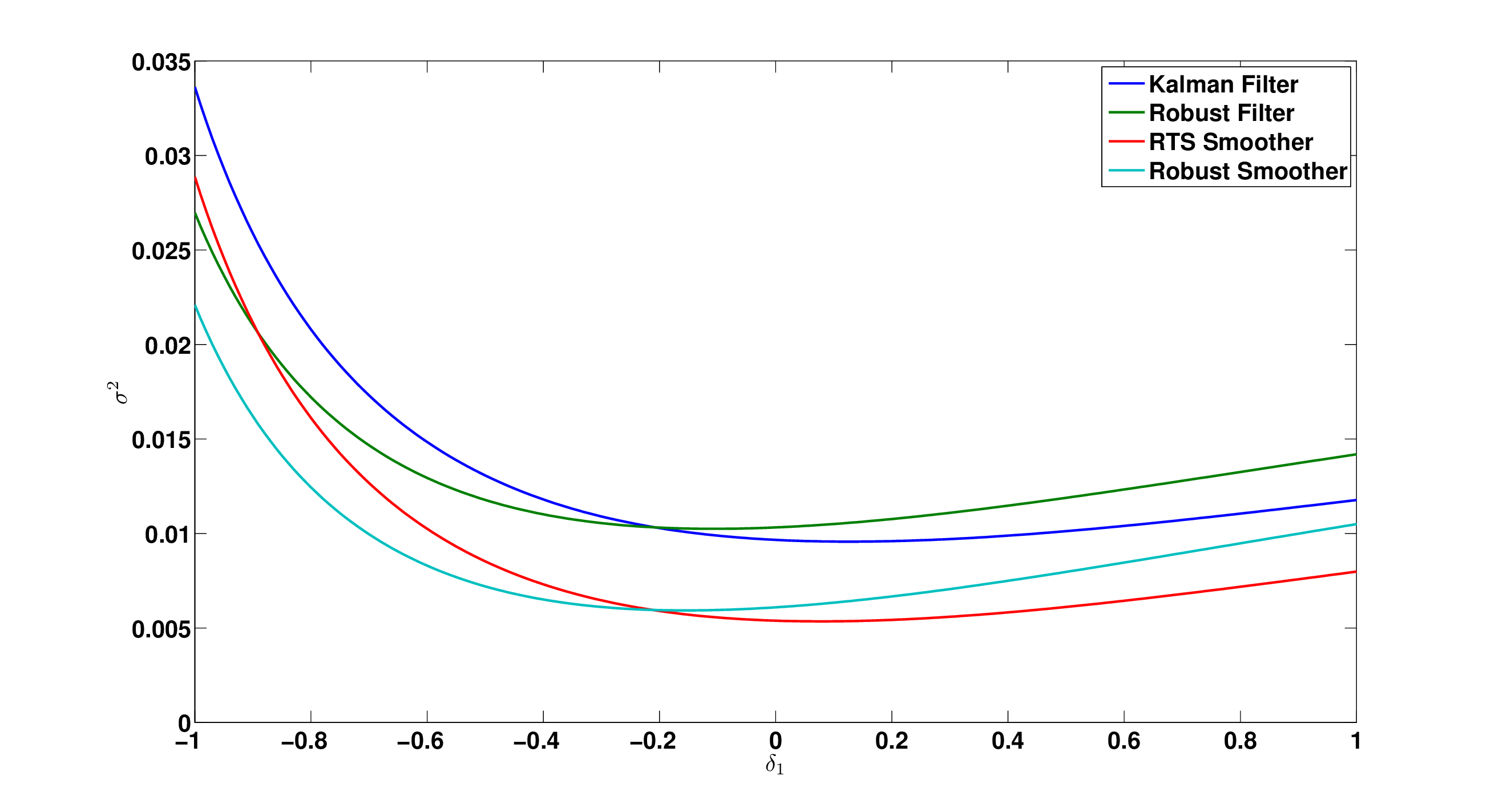}
\caption{Coherent State: Comparison of the smoothers for $\mu=0.7$.}
\label{fig:sql_mu_70}
\end{figure}

\begin{figure}[!htb]
\hspace*{-5mm}
\includegraphics[width=0.56\textwidth]{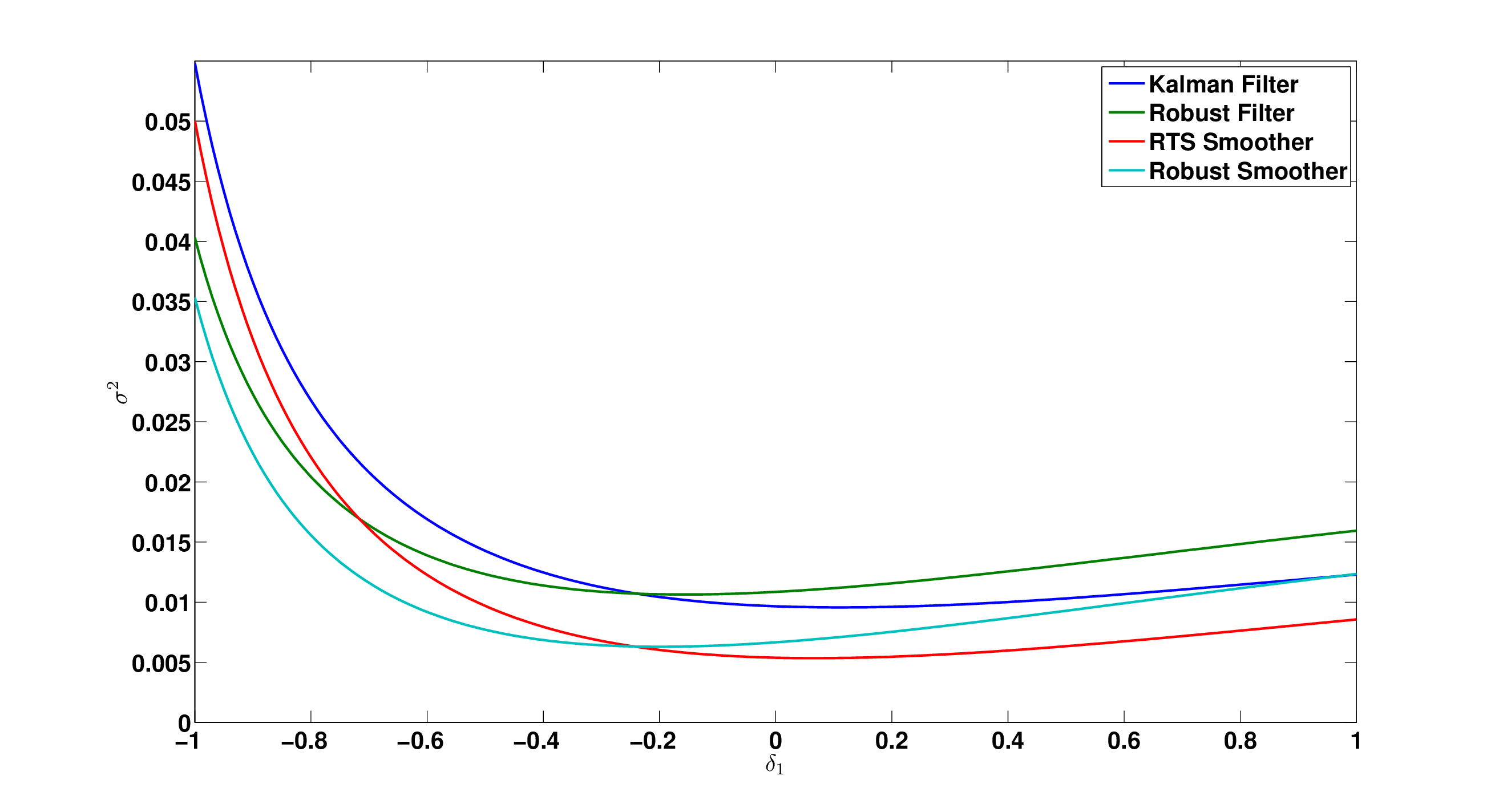}
\caption{Coherent State: Comparison of the smoothers for $\mu=0.8$.}
\label{fig:sql_mu_80}
\end{figure}

The estimation mean-square errors may be calculated, as described in Section \ref{sec:lyap_method} (i.e. $\Pi$), for the RTS smoother, and likewise for the robust smoother, as a function of the uncertain parameter $\delta$. The errors may similarly be computed (i.e. $\mathbf{\Pi_f}(1,1)$ in Section \ref{sec:lyap_method}) and plotted on the same graph for the forward Kalman filter alone and the forward robust filter alone, for comparison. Here, we use the nominal values of the parameters and choose different values for $\mu$. These values were used to generate plots of the errors versus $\delta$ to compare the performance of the robust smoother and the RTS smoother for the uncertain system. Figs. \ref{fig:sql_mu_50}, \ref{fig:sql_mu_70} and \ref{fig:sql_mu_80} show these plots for $\mu=0.5, 0.7$ and $0.8$, respectively.

Clearly, the RTS smoother behaves better than the robust smoother when $\delta=0$, as expected. However, in the worst-case scenario, the performance of the robust smoother is superior to that of the RTS smoother for all levels of $\mu$. Also, the robust smoother behaves better than the robust filter alone. Moreover, the improvement with the robust smoother over the optimal smoother is better with the resonant noise process considered here as compared to that with OU noise process considered in \cite{RPH2}. For example, while the worst-case improvement for $80\%$ uncertainty in the OU noise case was $\sim0.06$ dB, that in this resonant noise case is $\sim1.5$ dB.

\section{SQUEEZED STATE CASE}

While a coherent state has the same spread in both (amplitude and phase) quadratures, a squeezed state has reduced fluctuations in one of the two quadratures at the expense of increased fluctuations in the other. Here, we consider a phase-squeezed beam as in \cite{YNW}, whose phase is modulated with the resonant noise process. The beam is then measured by homodyne detection using a local oscillator, the phase of which is adjusted according to the filtered estimate $\phi_f(t)$.

The normalized homodyne output current $I(t)$ is given by
\begin{align}
I(t)dt &\simeq 2|\alpha|[\phi(t)-\phi_f(t)]dt + \sqrt{\overline{R}_{sq}}dW(t),\\
\overline{R}_{sq} &= \sigma_f^2 e^{2r_p}+(1-\sigma_f^2)e^{-2r_m},
\end{align}
where $|\alpha|$ is the amplitude of the input phase-squeezed beam, and $W(t)$ is a Wiener process arising from squeezed vacuum fluctuations. The parameter $\overline{R}_{sq}$ is determined by the degree of squeezing ($r_m \geq 0$) and anti-squeezing ($r_p \geq r_m$) and by $\sigma_f^2=\langle[\phi(t)-\phi_f(t)]^2\rangle$. The measurement model is \cite{RPH3}:
\begin{equation}
{\theta = \mathbf{Hx}+w,}
\end{equation}
where $\mathbf{H} = \left[\begin{array}{cc} 2|\alpha|/\sqrt{\overline{R}_{sq}} & 0 \end{array}\right]$.

Eqs. (\ref{eq:sys_model}) to (\ref{eq:robust_smoother_eqn}) are then modified accordingly.

\begin{figure}[!b]
\hspace*{-5mm}
\includegraphics[width=0.56\textwidth]{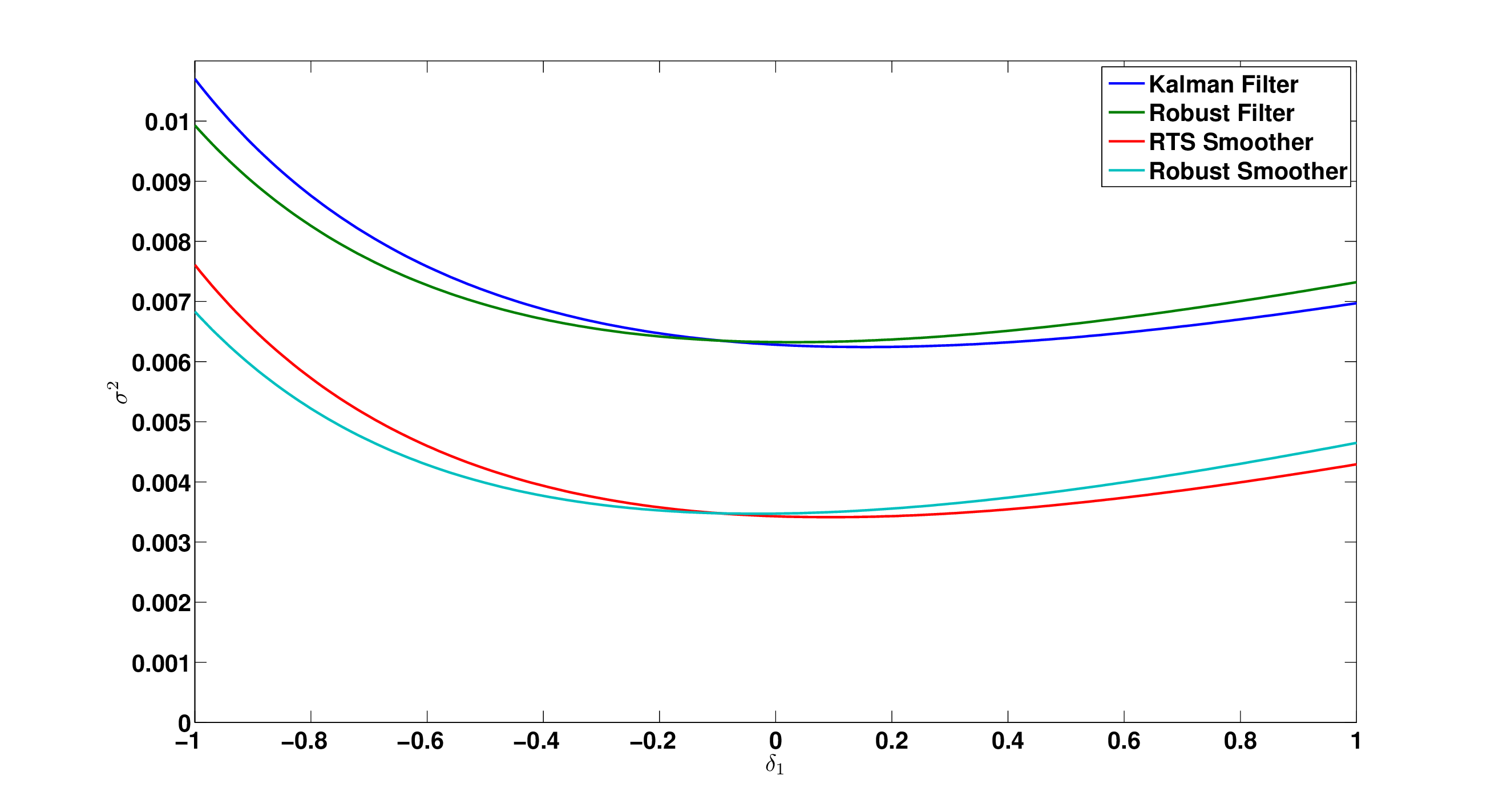}
\caption{Squeezed State: Comparison of the smoothers for $\mu=0.5$.}
\label{fig:csl_mu_50}
\end{figure}

\begin{figure}[!b]
\hspace*{-5mm}
\includegraphics[width=0.56\textwidth]{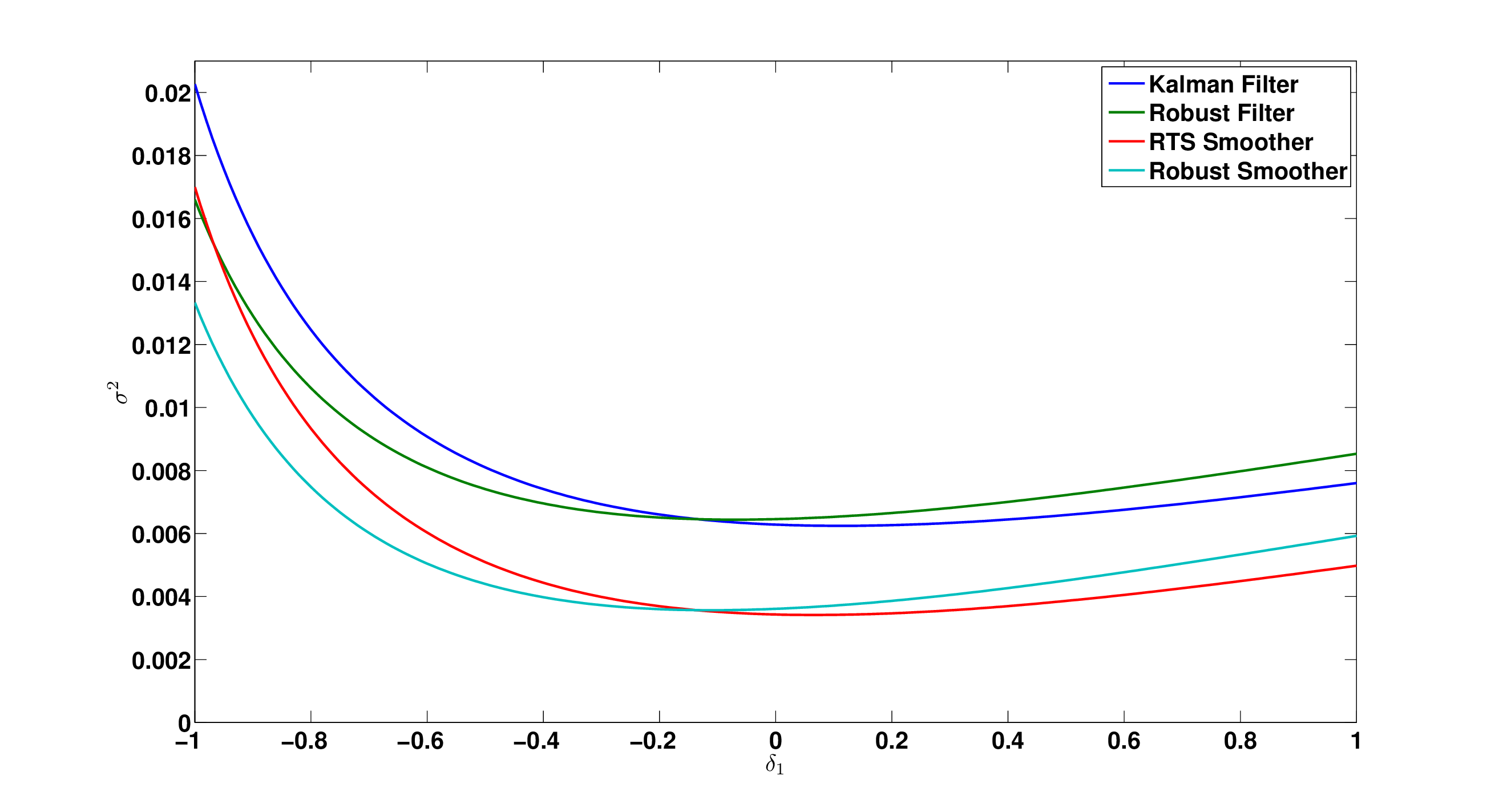}
\caption{Squeezed State: Comparison of the smoothers for $\mu=0.7$.}
\label{fig:csl_mu_70}
\end{figure}

\begin{figure}[!b]
\hspace*{-5mm}
\includegraphics[width=0.56\textwidth]{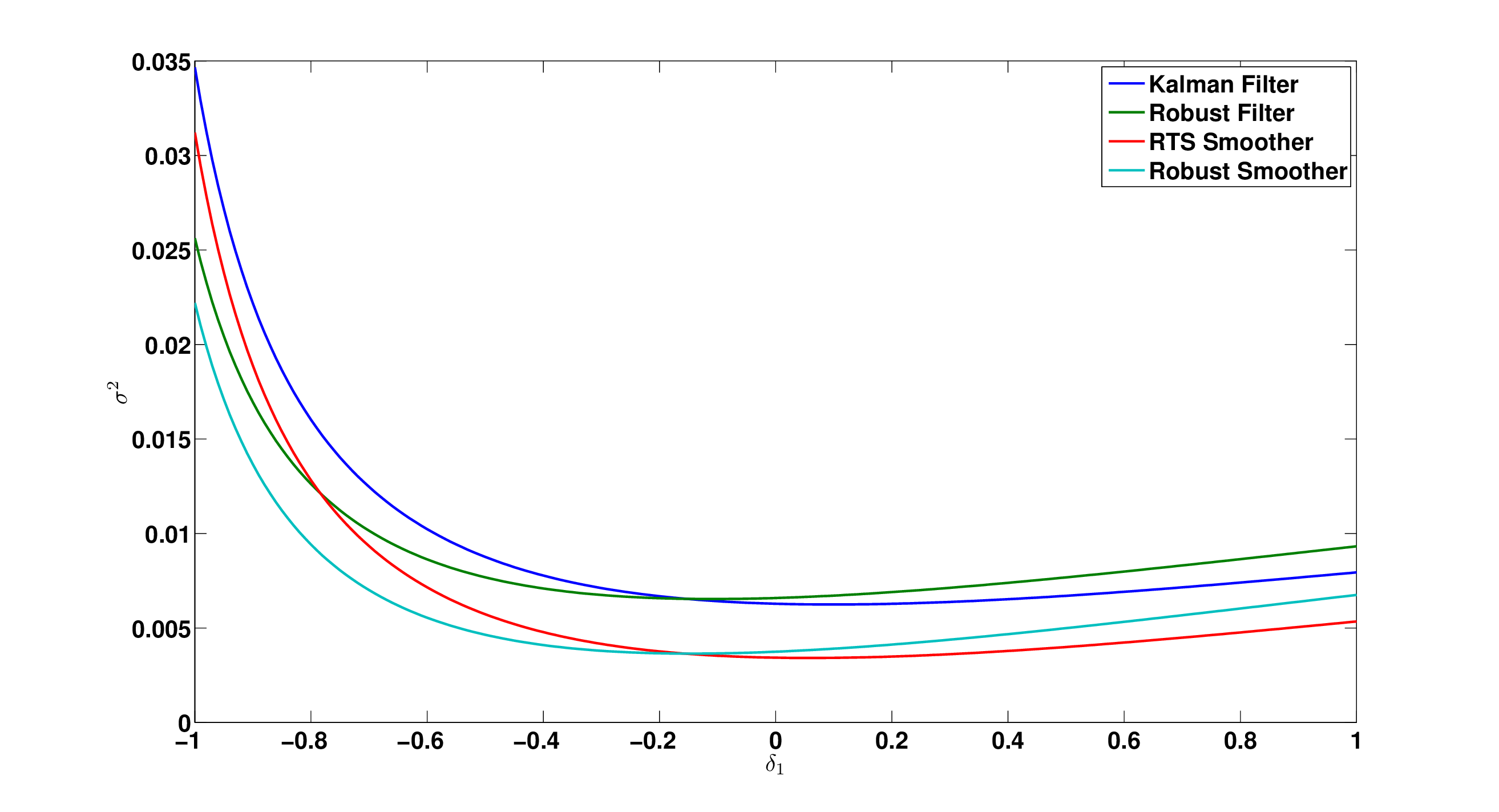}
\caption{Squeezed State: Comparison of the smoothers for $\mu=0.8$.}
\label{fig:csl_mu_80}
\end{figure}

We use the technique as described in Section \ref{sec:lyap_method} again to compute the estimation mean-square errors for the robust smoother and the RTS smoother, as a function of $\delta$. Again, the errors may as well be computed and plotted on the same graph for the forward Kalman filter alone and the forward robust filter alone, for comparison. Here, we use the nominal values of the parameters, $r_m=0.36$, $r_p=0.59$ and chosen values for $\mu$. These values were used to generate plots of the errors versus $\delta$ to compare the performance of the robust smoother and the RTS smoother for the uncertain system. Due to the implicit dependence of $\overline{R}_{sq}$ and $\sigma_f^2$, we compute the smoothed mean-square error in each case by running several iterations until $\sigma_f^2$ is obtained with an accuracy of $6$ decimal places. Figs. \ref{fig:csl_mu_50}, \ref{fig:csl_mu_70} and \ref{fig:csl_mu_80} show the plots for $\mu = 0.5, 0.7$ and $0.8$, respectively.

Note that the worst-case performance of the robust smoother is better than that of the RTS smoother for all levels of $\mu$. Also, the robust smoother behaves better than the robust filter alone. Moreover, the worst-case errors of the robust smoother is lower than those for the coherent state case considered earlier. For example, for $\mu=0.8$, there is a $\sim2$ dB improvement in the mean-square error of the robust smoother in the squeezed state case as compared to that in the coherent state case.

\section{CONCLUSION}

This work applies robust fixed-interval smoothing to homodyne phase estimation of coherent and squeezed states of light, when under the influence of a continuous-time resonant noise process. The robust smoother has been shown, as expected, to yield lower estimation errors than a robust filter alone in both the cases. More importantly, we have shown that for the uncertain system, the robust smoother performs better than the optimal smoother in the worst-case for both coherent and squeezed states, and the improvement so observed is better in this case of a resonant noise process than that observed (in earlier papers) in the case of an OU noise process. Also, the robust smoother provides superior accuracy in the estimate in the squeezed state case as compared to the coherent state case.

\bibliographystyle{IEEEtran}
\bibliography{IEEEabrv,bibliography}

\end{document}